\documentclass[12pt]{article}
\usepackage {amsthm,amsmath, amssymb}
\def\ep{\varepsilon}
\begin{document}
\centerline {\bf\Large On a set of transformations}

\centerline {\bf\Large of Gaussian random functions}
\bigskip

\centerline {A.I.~Nazarov\footnote{Partially supported
by grant NSh.227.2008.1 and by RFBR grant 07-01-00159.}}

\centerline {Dept. of Mathematics and Mechanics,}
\centerline {St.Petersburg State University }

\centerline {\small e-mail:\ an@AN4751.spb.edu}

{\it\small We consider a set of one-dimensional transformations of Gaussian
random functions. Under natural assumptions we obtain a connection between
$L_2$-small ball asymptotics of the transformed function and of the original
one. Also the explicit Karhunen -- Lo\'eve expansion is obtained for a proper
class of Gaussian processes}.



\section{Introduction and Main Lemma}

Recently P.~Deheuvels \cite{De} showed that for a standard Brownian bridge
$B(t)$, $0\leq t\leq 1$, the distributional equality
$${\cal Y}_K(t)\,\stackrel {d}{=}\,{\cal Y}_{2-K}(t), \qquad t\in [0,1], \
K\in\mathbb R,\eqno(1.1)$$
holds true (here ${\cal Y}_K(t)=B(t)-6Kt(1-t)\int_0^1B(s)\,ds$). Moreover, he
obtained the explicit Karhunen -- Lo\'eve (KL) expansion for the process
${\cal Y}_1(t)$.

We introduce one-parameter sets of transformations for zero mean-value Gaussian
random functions. These transformations satisfy a relation generalizing (1.1).
If the $L_2$-norm of original function is finite a.s. we derive the explicit
relation between exact asymptotics of $L_2$-small ball probabilities for the
transformed function and for the original one. For one-dimensional processes
generating boundary value problems to ordinary differential equations we
obtain also the KL-expansion for transformed process provided the KL-expansion
for original process is known.\medskip

Let us consider a zero mean-value Gaussian random function $X(x)$, $x\in\overline{\cal O}$,
with the covariance $G_X(x,y)={\mathbb E}X(x)X(y)$, $x,y\in\overline{\cal O}$. For simplicity we suppose that ${\cal O}$ is a bounded domain in $\mathbb R^n$.

Let $\varphi$ be a locally summable function in ${\cal O}$. Suppose that the function
$$\psi(x)=\int\limits_{\cal O} G_X(x,y)\varphi(y)\,dy\eqno(1.2)$$
is well defined a.e. in ${\cal O}$, $\psi\not\equiv0$, and
$$q=\int\limits_{\cal O}\psi(u)\varphi(u)\,du=\int\limits_{\cal O}
\int\limits_{\cal O} G_X(u,v)\varphi(u)\varphi(v)\,dudv<\infty.\eqno(1.3)$$

We define a set of Gaussian functions
$${\cal X}_{\varphi,\alpha}(x)=X(x)-
\alpha\psi(x)\int\limits_{\cal O} X(u)\varphi(u)\,du,\quad
x\in\overline{\cal O}.\eqno(1.4)$$

{\bf MAIN LEMMA}. The covariance of ${\cal X}_{\varphi,\alpha}$ is
$${\cal G}_{\varphi,\alpha}(x,y)=G_X(x,y)+
Q\psi(x)\psi(y),\eqno(1.5)$$
where $Q=q\alpha^2-2\alpha$.\medskip

{\bf Proof}. The formula (1.5) can be checked by direct computation due to (1.2).
\hfill$\square$\medskip

{\bf Corollary 1}. The equality
$${\cal X}_{\varphi,\alpha}(x)\,\stackrel {d}{=}\,
{\cal X}_{\varphi,\frac 2q-\alpha}(x),\quad x\in\overline{\cal O}.$$
holds for the processes (1.4). In particular, ${\cal X}_{\varphi,\frac 2q}(x)\,\stackrel {d}{=}\,X(x)$.
\medskip

{\bf Corollary 2}. Let $\widehat\alpha=\frac 1q$. Then

1. The following identity holds true a.s.:
$$\int\limits_{\cal O}{\cal X}_{\varphi,\widehat\alpha}(x)\varphi(x)\,dx=0.$$

2. The process ${\cal X}_{\varphi,\widehat\alpha}(x)$ and the r.v.
$\int_{\cal O} X(u)\varphi(u)du$ are independent.

3. If $\varphi\in L_2({\cal O})$ then the integral operator with the kernel function
${\cal G}_{\varphi,\widehat\alpha}(x,y)$ has a zero eigenvalue with the eigenfunction $\varphi$.\medskip

{\bf Proof}. All three statements follow from relations
$$\int\limits_{\cal O}{\cal X}_{\varphi,\alpha}(x)\varphi(x)\,dx=
\int\limits_{\cal O} X(u)\varphi(u)\,du\cdot(1-q\alpha);$$
$${\mathbb E}{\cal X}_{\varphi,\alpha}(x)
\int\limits_{\cal O} X(u)\varphi(u)\,du=\psi(x)\cdot(1-q\alpha),
$$
which can be easily checked.\hfill$\square$\medskip

{\bf Remark}. Trivially the process ${\cal Y}_K$ coincides with ${\cal B}_{\varphi, \alpha}$ for $\varphi\equiv1$, $\alpha=12K$. Hence, Lemma 2.2 and Corollaries 2.1, 2.2 from \cite{De} are particular cases of our statements.

\section{The asymptotics of small ball probabilities in $L_2$}

Now we suppose that
$$\|X\|_2^2\equiv\int\limits_{{\cal O}}X^2(x)\,dx<\infty\qquad\mbox{ a.s.}\eqno(2.1)$$
Then the process $X$ admits the KL-expansion
$$X(x)\stackrel {d}{=}\,\sum\limits_{k=1}^{\infty}\sqrt
{\lambda_k}u_k(x)\xi_k,\qquad x\in{\cal O},\eqno(2.2)$$
where $\xi_k$ is a sequence of independent standard Gaussian r.v. while $\lambda_k>0$ and $u_k$ are, respectively, eigenvalues and (normalized in $L_2({\cal O})$) eigenfunctions of the integral operator ${\mathfrak G}$ with the kernel function $G_X(x,y)$. Moreover, (2.1) implies $\sum_k\lambda_k<\infty$ i.e the operator ${\mathfrak G}$ belongs to the kernel class ${\mathfrak S}_1$. Remark that the series (2.2) converges in $L_2({\cal O})$ a.s.

The relation (2.2) implies
$$\|X\|_2^2\stackrel {d}{=}\,\sum\limits_{k=1}^{\infty}\lambda_k\xi_k^2.$$
Therefore, having the eigenvalues $\lambda_k$ in hands one can obtain some information on the distribution of $\|X\|_2^2$. In particular, one can derive the exact asymptotics of small ball probabilities in $L_2$ i.e. describe the behavior of the probability ${\mathbb P}\{\|X\|_2\leq \ep\}$ as $\ep\to0$.

In \cite{Syt}, a solution of the small ball behavior problem was obtained in abstract form.
Then  many authors simplified formulas for small ball probabilities under various assumptions;
see \cite{DLL} and references therein.

In papers \cite{NN}, \cite{Na}, \cite{Na1} a new approach was delivered. This approach provides the exact $L_2$-small ball asymptotics for a zero mean-value Gaussian process if its covariance is
the Green function for an ordinary differential operator. In more general case the problem cannot be
solved completely yet. However, in the case under consideration the transformed operator is a one-dimensional perturbation  of the original one. So, we can derive the small ball asymptotics for
the process ${\cal X}_{\varphi,\alpha}$ in terms of the small ball asymptotics for the original
process.\medskip

{\bf Remark}. If the operator ${\mathfrak G}$ has a nontrivial null-space $U_0$ then it is
contained in the null-space of the operator with the kernel function
${\cal G}_{\varphi,\alpha}(x,y)$ due to obvious relation $\psi\perp U_0$. Therefore, we can apply all
the arguments in an orthogonal complement to $U_0$ (i.e. in the image of ${\mathfrak G}$). So, without
loss of generality one can assume $\{u_k\}$ in (2.2) being an orthogonal basis.\medskip

{\bf Theorem 1}. Let the process $X$ satisfy (2.1). Suppose that a function
$\varphi\in L_{1,loc}({\cal O})$ satisfies (1.3). If $\alpha\ne \frac 1q$ then, as $\ep\to0$,
$${\mathbb P}\{\|{\cal X}_{{\bf \varphi}, \alpha}\|_2\leq\ep\} \sim
 \frac 1{|1-\alpha q|} \cdot
 {\mathbb P}\{\|X\|_2\leq\ep\}.\eqno(2.3)
$$

{\bf Proof}. By comparison theorem (\cite{Li}; see also \cite{GHT}),
$${\mathbb P}\{\|{\cal X}_{{\bf \varphi}, \alpha}\|_2\leq\ep\} \sim
{\mathbb P}\{\|X\|_2\leq\ep\}\cdot \Bigl(\prod_{k=1}^{\infty}
\frac {\lambda_k}{\widetilde\lambda_k}\Bigr)^{1/2}, \eqno(2.4)$$
where $\widetilde\lambda_k$ are the eigenvalues of the integral operator with kernel function
${\cal G}_{\varphi,\alpha}(x,y)$. Remark that due to the minimax principle,
see, e.g., \cite[\S10.2]{BS}) the sequences $\lambda_k$ and $\widetilde\lambda_k$ interlace.
In particular, this implies the convergence of the series $\sum_k \widetilde\lambda_k$.\medskip

By definition, put $\mu_k=\lambda_k^{-1}$, $\widetilde\mu_k=\widetilde\lambda_k^{-1}$.
Consider the Fredholm determinants for the kernels $G_X$ and ${\cal G}_{\varphi,\alpha}$:
$${\cal F}(z)=\prod_{k=1}^{\infty}\left(1-\frac {z}{\mu_k}\right); \qquad
\widetilde {\cal F}(z)=\prod_{k=1}^{\infty}\left(1-\frac
{z}{\widetilde\mu_k}\right).$$
Since the series
$\sum\limits_k\mu_k^{-1}$ and $\sum\limits_k \widetilde\mu_k^{-1}$
converge, these canonical Hadamard's products converge for all
$z\in\mathbb C$. By (1.5) the following relation holds
true\footnote {Note that (2.5) is a particular case of the
transformation formula for the Fredholm determinant under
finite-dimensional perturbation of the operator. In the literature
on statistics this formula usually is attributed to \cite{Su}.
Statistical applications in \cite{Su} seem to be new while the
formula itself was really obtained in \cite{Ba} and was well known
as in computational methods as in spectral theory, even in more
general situation, see, e.g., \cite[Ch.II, 4.6]{KK} and
\cite[Sec.106]{AG}.}:
$$\widetilde {\cal F}(z)={\cal F}(z)\cdot\Bigl(1+Q\sum_{k=1}^{\infty}
\frac {a_k^2\mu_k}{1-\frac {\mu_k}{z}}\Bigr),\eqno(2.5)$$
where $a_k$ are the Fourier coefficients of the function $\psi$ with respect to the system $\{u_k\}$.

The Jensen theorem, see \cite[\S3.6]{Ti}, provides
$$\prod_{k=1}^{\infty}\frac {\mu_k}{\widetilde\mu_k}=\lim\limits_{|z|\to\infty}
\exp\left(\frac 1{2\pi}\int\limits_0^{2\pi}\ln\left(\frac
{|\widetilde {\cal F}(z)|}{|{\cal F}(z)|}\right)d\arg(z)\right).\eqno(2.6)$$
Formula (2.5) and Lemma 5.1 show that the last limit equals $\big|1+Q\sum_k a_k^2\mu_k\big|$.
But $\psi=\sum_ka_ku_k$ implies $\varphi=\sum_k\mu_ka_ku_k$ and, therefore,
$$\sum_{k=1}^{\infty}a_k^2\mu_k=\int\limits_{\cal O}\psi(u)\varphi(u)\,du=q.
\eqno(2.7)$$
Substituting (2.6) into (2.4) we obtain (2.3).\hfill$\square$\medskip

Now we consider the critical case $\widehat\alpha=\frac 1q$.\medskip

{\bf Theorem 2}. Let the process $X$ satisfy (2.1), and let $\widehat\alpha=\frac 1q$.
If $\varphi\in L_2({\cal O})$ then, as $\ep\to0$,
$${\mathbb P}\{\|{\cal X}_{{\bf \varphi}, \widehat\alpha}\|_2\leq\ep\}
\sim \frac {\sqrt {q}}{\|\varphi\|_2} \cdot \sqrt {\frac 2{\pi}}\cdot
\int\limits_0^{\ep^2}\frac d{dt}{\mathbb P}\{\|X\|_2\leq t\}\ \frac {dt}
{\sqrt {\ep^2-t^2}}.\eqno(2.8)$$

{\bf Proof}. We introduce three distribution functions:
$$\gathered
F(r)={\mathbb P}\{\sum\limits_{k=1}^{\infty}\lambda_k\xi_k^2\leq r\}
={\mathbb P}\{\| X \|_2\leq \sqrt{r}\}; \\
\widetilde F(r)=
{\mathbb P}\{\sum\limits_{k=1}^{\infty}\widetilde\lambda_k\xi_k^2\leq r\}
={\mathbb P}\{\|{\cal X}_{{\bf \varphi}, \widehat\alpha}\|_2\leq \sqrt{r}\}; \\
F_1(r)={\mathbb P}\{\sum\limits_{k=2}^{\infty}\lambda_k\xi_k^2\leq r\}.
\phantom{={\mathbb P}\{\| X \|_2\leq \sqrt{r}\}}
\endgathered$$

Similarly to the previous theorem we have, as $r\to0$,
$$\widetilde F(r)\sim F_1(r)\cdot \Bigl(\prod_{k=2}^{\infty}
\frac {\widetilde\mu_k}{\mu_{k-1}}\Bigr)^{1/2}.$$
The Jensen theorem provides
$$\prod_{k=2}^{\infty}\frac {\mu_{k-1}}{\widetilde\mu_k}=
\lim\limits_{|z|\to\infty}\exp\left(\frac 1{2\pi}\int\limits_0^{2\pi}
\ln\left|\Bigl(1-\frac z{\mu_1}\Bigr)\cdot\frac
{\widetilde {\cal F}(z)}{{\cal F}(z)}\right|d\arg(z)\right).\eqno(2.9)$$
The assumption $\widehat\alpha=\frac 1q$ implies $Q=-\frac 1q$. Hence, due to (2.5) and (2.7) the
expression under log sign can be rewritten as follows:
$$\left|\Bigl(1-\frac z{\mu_1}\Bigr)\cdot\Bigl(1+Q\sum_{k=1}^{\infty}
\frac {a_k^2\mu_k}{1-\frac {\mu_k}{z}}\Bigr)\right|=\left|\frac {\frac 1{\mu_1}-\frac 1z}
q\cdot \sum_{k=1}^{\infty}\frac {a_k^2\mu_k^2}{1-\frac {\mu_k}{z}}\right|.$$
By Lemma 5.1, the limit in (2.9) equals
$$\frac 1{\mu_1q}\sum_{k=1}^{\infty}a_k^2\mu_k^2=\frac {\|\varphi\|_2^2}{\mu_1q}.$$
This gives
$$\widetilde F(r)\sim F_1(r)\cdot\frac {\sqrt {q\mu_1}}{\|\varphi\|_2}.
\eqno(2.10)$$

Further, obviously, $F(r)=(F_1*f)(r)$, where
$$f(x)=\frac d{dx}{\mathbb P}\{\lambda_1\xi^2\leq x\}=\frac {\exp\bigl(-\frac
x{2\lambda_1}\bigr)}{\sqrt{2\pi\lambda_1x}}.$$
By the Laplace transform we obtain a solution of this convolution equation:
$$F_1(r)=\sqrt {\frac {2\lambda_1}{\pi}}\exp\bigl(-\frac r{2\lambda_1}\bigr)
\int\limits_0^r\Bigl(F(x)\exp\bigl(\frac x{2\lambda_1}\bigr)\Bigr)'\frac
{dx}{\sqrt{r-x}}.\eqno(2.11)$$

{\bf Lemma 2.1}. $F(x)=o(F'(x))$ as $x\to+0$.\medskip

{\bf Proof}. Remark that $F'$ is absolutely continuous on $\mathbb R$, if the sum (2.2) contains
at least three nonzero summands. We claim that $F''>0$ in a right half-neighborhood of the origin.
Indeed, this property can be directly checked for three summands and is easily conserved under adding
a new summand; moreover, the radius of a neighborhood is not decreasing. Consequently, this property
is conserved also for the infinite sum.

Since $F$ is convex in a right half-neighborhood of the origin, $F'(x)\ge F(x)/x$ in this
half-neighborhood.\hfill$\square$\medskip

We continue the proof of Theorem 2. By Lemma 2.1,  we obtain from (2.11)
$$F_1(r)\sim\sqrt {\frac {2\lambda_1}{\pi}}\int\limits_0^r F'(x)\cdot\frac
{\exp\bigl(-\frac {r-x}{2\lambda_1}\bigr)dx}{\sqrt{r-x}}\sim\sqrt
{\frac {2\lambda_1}{\pi}}\int\limits_0^r F'(x)\ \frac {dx}{\sqrt{r-x}}.
\eqno(2.12)$$
Substituting  (2.12) into (2.10) and changing the variable $r=\ep^2$ we arrive at (2.8).
\hfill$\square$\medskip

{\bf Remark}. By misuse of language, we can interpret (2.10) as
follows. For $\alpha\ne\frac 1q$ the sequences $\mu_k$ and
$\widetilde\mu_k$ have the same asymptotics. Under the conditions
of Theorem 2, the statement 3 of Corollary 2 provides vanishing of
an eigenvalue $\widetilde\lambda_k$. This generates a confusion in
the enumeration of $\widetilde\mu_k$; removing of a summand from
(2.2) we reestablish the corresponding of enumerations.\medskip

\section{Karhunen -- Lo\'eve expansion}

Now we suppose that $n=1$, ${\cal O}=(0,a)$ is an interval, and
the covariance $G_X(t,s)$, $t,s\in[0,a]$, is the Green function of
a self-adjoint operator $L_X$ in the space $L_2(0,a)$, generated
by differential expression of order $2\ell$
$$L_Xu\equiv
(-1)^{\ell}u^{({2\ell})}+
\left(p_{\ell-1}u^{({\ell}-1)}\right)^{({\ell}-1)}+\dots+p_0u,\eqno(3.1)$$
and $2\ell$ boundary conditions. We recall that by definition
$G_X$ for any $s\in(0,a)$ satisfies the equation
$L_XG_X=\delta(t-s)$ in the sense of distributions, and satisfies
boundary conditions. Without loss of generality we assume $a=1$.

By ${\cal D}(L_X)$ denote the image of integral operator with the kernel function $G_X(t,s)$.
Then it is easy to see that the inverse operator is just $L_X$ with the domain ${\cal D}(L_X)$.
In particular, if $\varphi\in L_2(0,1)$ then $\psi\in{\cal D}(L_X)$, and $L_X\psi=\varphi$.

Assume for simplicity that $p_j\in{\cal C}^j[0,1]$. Then ${\cal D}(L_X)$ coincides with the set
of functions which belong to $W^{\ell}_2(0,1)$ and satisfy boundary conditions. By (1.5) we
obtain that the covariance of transformed process ${\cal X}_{\varphi,\alpha}$ satisfies the equation
$$L_X{\cal G}_{\varphi,\alpha}=\delta(t-s)+Q\varphi(t)\psi(s)\eqno(3.2)$$
(recall that $Q=q\alpha^2-2\alpha$) and satisfies the same boundary conditions.

Suppose we know KL-expansion (2.2) for the original process. Then, obviously, $u_k$ are the
eigenfunctions and $\mu_k=\lambda_k^{-1}$ are the eigenvalues of the boundary value problem
$$L_Xu=\mu u, \qquad u\in{\cal D}(L_X).\eqno(3.3)$$
In practice, the eigenfunctions of the problem (3.3) can be found analytically only if we know
the fundamental system of solutions to the equation $L_Xv-\mu v=0$ for arbitrary $\mu\in\mathbb R$.
We take advantage of this to give an algorithm of derivation of KL-expansion for the process
${\cal X}_{\varphi,\alpha}$. Note that in a particular case this algorithm was used in \cite{KKW}
(see below the example {\bf5}).

From (3.2) we obtain the boundary value problem for eigenfunctions of integral operator
with the kernel function ${\cal G}_{\varphi,\alpha}(t,s)$:
$$L_Xu=\mu u+\mu Q\varphi\int\limits_0^1u(s)\psi(s)\,ds,
\qquad u\in{\cal D}(L_X).\eqno(3.4)$$

Let $\mu$ be an unknown parameter. By the Lagrange method we can construct a particular
solution to the equation $L_X\eta-\mu\eta=\varphi$. Then a general solution to the equation (3.4)
can be written as follows:
$$u=c_0\eta+c_1v_1+c_2v_2+\dots+c_{2\ell}v_{2\ell},$$
where $v_1,\dots,v_{2\ell}$ form a fundamental system for the equation $L_Xv-\mu v=0$. Substituting
$u$ into the boundary conditions we obtain $2\ell$ equations for the constants
$c_0,c_1,\dots,c_{2\ell}$. One more equation follows from the equality of the coefficients at
$\varphi$ in (3.4):
$$\frac {c_0}{\mu Q}= c_0\int\limits_0^1\eta(s)\psi(s)\,ds+
c_1\int\limits_0^1v_1(s)\psi(s)\,ds+\dots+
c_{2\ell}\int\limits_0^1v_{2\ell}(s)\psi(s)\,ds.\eqno(3.5)$$
Eigenvalues of the problem (3.4) are roots of the determinant of
the obtained homogenous system while eigenfunctions are its
nontrivial solutions.\medskip

Let us show some examples.\medskip

{\bf Example 1}. Let $X=W$ be a standard Wiener process, and let $\varphi\equiv1$. Then
$$\psi(s)=\int\limits_0^1\min(t,s)\,ds=\frac {2t-t^2}2;\qquad
q=\int\limits_0^1\frac {2t-t^2}2\,dt=\frac 13.$$
The relation (3.4) reads as follows:
$$-u''=\mu u +\frac {\mu Q}2\int\limits_0^1u(s)\,(2s-s^2)\,ds,\qquad
u(0)=u'(1)=0.$$
A general solution of this equation is
$$u(t)=c_0+c_1\cos(\omega t)+c_2\sin(\omega t),\qquad
\omega=\mu^{1/2}.\eqno(3.6)$$
Substituting it into the boundary conditions and into (3.5) we derive the equation for eigenvalues:
$$Q\sin(\omega)=\cos(\omega)\cdot(Q\omega+(1+Q/3)\omega^3).\eqno(3.7)$$
For $Q=0$ we obtain $\cos(\omega)=0$. This is natural result because in
this case we have a conventional Wiener process.

We remark also that $1+Q/3=(1-\alpha/3)^2\ge0$. For $\alpha=3$ the equation (3.7) is reduced to
$\tan(\omega)=\omega$.

Let $\omega_k$ be positive roots of (3.7) enumerated in the increasing order.
Then we put in KL-expansion of the process ${\cal W}_{{\bf 1}, \alpha}$
$$\widetilde\lambda_k=\omega_k^{-2},\qquad
\widetilde u_k(t)=\gamma_k(\cos(\omega_kt)-1+\tan(\omega_k)\sin(\omega_k t)),$$
where $\gamma_k$ are the normalizing constants.\medskip

{\bf Example 2}. Let $X=B$ be a standard Brownian bridge, and let $\varphi\equiv1$. Then
$\psi(t)=\frac {t-t^2}2$, $q=\frac 1{12}$.
The relation (3.4) reads as follows:
$$-u''=\mu u +\frac {\mu Q}2\int\limits_0^1u(s)\,(s-s^2)\,ds,\qquad
u(0)=u(1)=0.\eqno(3.8)$$
A general solution of this equation is given by (3.6). Substituting it into the boundary conditions and into (3.5) we derive the equation for eigenvalues:
$$\left[
\begin{array}{rcl}
\sin(\tau)&=&0;\\
Q\sin(\tau)&=&\cos(\tau)\cdot(Q\tau+(4+Q/3)\tau^3),
\end{array}\right.
\qquad \tau=\frac\omega 2.
\eqno(3.9)$$
For $Q=0$ two equations in (3.9) can be merged into $\sin(\omega)=0$, as is expected.
We remark also that $4+Q/3=(2-\alpha/6)^2\ge0$. For $\alpha=12$ the second equation in
(3.9) is reduced to $\tan(\tau)=\tau$. This case was considered in \cite[Theorem 1.2]{De}.

A sequence
$$Q_k=-\frac {12}{1+\frac 3{(k\pi)^2}}\searrow -12$$
(any $Q_k$ corresponds to two values of $\alpha$) has a curious property. Namely, for $Q=Q_k$ the
$k$-th roots of both equations in (3.9) coincide. Thus, for $Q=Q_1$ the least eigenvalue of the
problem (3.8) is multiple. This effect is impossible for a conventional Sturm -- Liouville problem.

Let $\tau_k$ be positive roots of the second equation in (3.9) enumerated in the increasing order.
Then we put in KL-expansion of the process ${\cal B}_{{\bf 1}, \alpha}$
$$\begin{array}{ll}
\widetilde\lambda_k=(2k\pi)^{-2},\qquad& \widetilde u_k(t)=\sqrt{2}\sin(2\pi kt);\\
\widetilde{\widetilde\lambda}_k=(2\tau_k)^{-2},\qquad &
\widetilde{\widetilde u}_k(t)=\gamma_k(\cos(2\tau_kt)-1+\tan(\tau_k)\sin(2\tau_k t)),
\end{array}$$
where $\gamma_k$ are the normalizing constants.\medskip

{\bf Example 3}. Let $X=B$, and let $\varphi(t)=t(1-t)$. Then
$\psi(t)=\frac {t-2t^3+t^4}{12}$, $q=\frac {17}{5040}$.
The relation (3.4) reads as follows:
$$-u''=\mu u +\frac {\mu Q}{12}t(1-t)\int\limits_0^1u(s)\,(s-2s^3+s^4)\,ds,
\qquad u(0)=u(1)=0.$$
A general solution of this equation is
$$u(t)=c_0(t-t^2+{\textstyle\frac 2{\omega^2}})+
c_1\cos(\omega t)+c_2\sin(\omega t),\qquad
\omega=\mu^{1/2}.
$$
Substituting it into the boundary conditions and into (3.5) we derive the equation for eigenvalues:
$$\left[
\begin{array}{rcl}
\sin(\tau)&=&0;\\
Q\sin(\tau)&=&\cos(\tau)\cdot(Q\tau+Q\tau^3/3+2Q\tau^5/15+16(1+Qq)\tau^7),
\end{array}\right.
\qquad \tau=\frac\omega 2.
\eqno(3.10)$$
For $Q=0$ two equations in (3.10) can be merged into $\sin(\omega)=0$.
We remark also that $1+Qq=(1-q\alpha)^2\ge0$.

Let $\tau_k$ be positive roots of the second equation in (3.10) enumerated in the increasing order.
Then we put in KL-expansion of the process ${\cal B}_{\varphi, \alpha}$
$$\begin{array}{ll}
\widetilde\lambda_k=(2k\pi)^{-2},\qquad&\widetilde u_k(t)=\sqrt{2}\sin(2\pi kt);\\
\widetilde{\widetilde\lambda}_k=(2\tau_k)^{-2},\qquad &
\widetilde{\widetilde u}_k(t)=\gamma_k(\cos(2\tau_kt)-1+ 2\tau_k^2(t^2-t)+
\tan(\tau_k)\sin(2\tau_k t)),
\end{array}$$
where $\gamma_k$ are the normalizing constants.\medskip

{\bf Example 4}. Let
$$X(t)={\overline W}_1(t)=\int\limits_0^t
\Bigl(W(s)-\int\limits_0^1W(u)\,du\Bigr)\,ds$$
be an integrated centered Wiener process. Its covariation $G_{{\overline W}_1}$ is
the Green function of the operator $L_{{\overline W}_1}=L_B^2$, see, e.g., \cite{HN} and
 \cite[Prop. 5.4]{NN}.

If $\varphi\equiv1$, then $\psi(t)=\frac {t-2t^3+t^4}{24}$, $q=\frac 1{120}$,
and the relation (3.4) reads as follows:
$$u^{\scriptscriptstyle IV}=\mu u +\frac {\mu Q}{24}
\int\limits_0^1u(s)\,(s-2s^3+s^4)\,ds,
\qquad u(0)=u(1)=u''(0)=u''(1)=0.$$
A general solution of this equation is
$$u(t)=c_0+c_1\cos(\omega t)+c_2\sin(\omega t)+c_3\cosh(\omega t)+
c_4\sinh(\omega t),\qquad \omega=\mu^{1/4}.
$$
Substituting it into the boundary conditions and into (3.5) we derive the equation for eigenvalues:
$$\left[
\begin{array}{rcl}
\sin(\tau)&=&0;\\
Q(\sin(\tau)+\cos(\tau)\tanh(\tau))&=&
\cos(\tau)\cdot(2Q\tau+(32+4Q/15)\tau^5),
\end{array}\right.
\qquad \tau=\frac\omega 2.
\eqno(3.11)$$
For $Q=0$ two equations in (3.10) can be merged into $\sin(\omega)=0$.
We remark also that $32+4Q/15=32(1-q\alpha)^2\ge0$.

Let $\tau_k$ be positive roots of the second equation in (3.11) enumerated in the increasing order.
Then we put in KL-expansion of the process ${\cal X}_{{\bf 1}, \alpha}$
$$\begin{array}{ll}
\widetilde\lambda_k=(2k\pi)^{-4},\qquad&\widetilde u_k(t)=\sqrt{2}\sin(2\pi kt);\\
\widetilde{\widetilde\lambda}_k=(2\tau_k)^{-4},\qquad &
\widetilde{\widetilde u}_k(t)=\gamma_k(\cos(2\tau_kt)+\cosh(2\tau_kt)-2+\\
& +\tan(\tau_k)\sin(2\tau_k t) -\tanh(\tau_k)\sinh(2\tau_k t)),
\end{array}$$
where $\gamma_k$ are the normalizing constants.\medskip

It should be noted that the roots of the first equation in (3.9)-(3.11) are independent of $Q$.
The reason is that corresponding eigenfunctions $\sin(2\pi kt)$, $k\in\mathbb N$,
are orthogonal to $\psi$ in $L_2(0,1)$, and the last term in (3.4) vanishes.\medskip

{\bf Example 5}. Let $X=B$, and let
$$\varphi(t)=\frac 1{\phi(\mbox{\boldmath $\Phi$}^{-1}(t))},\qquad
\mbox{ЦДЕ}\qquad \phi(t)=\frac 1{\sqrt{2\pi}}\exp(-t^2/2);\qquad
\mbox {\boldmath $\Phi$}(x)=\int\limits_{-\infty}^x\phi(t)\,dt.$$
This example is important in statistics, see \cite{KKW}, \cite{Su}. It is evident that
$\varphi\not\in L_2(0,1)$. However, the direct calculation gives
$\psi=\phi(\mbox{\boldmath $\Phi$}^{-1})=\frac 1{\varphi}$, $q=1$. Hence all the statements
of \S1 are applicable here.

Further, the relation (3.4) reads as follows:
$$-u''=\mu u +\frac {\mu Q}{\phi(\mbox{\boldmath $\Phi$}^{-1}(t))}
\int\limits_0^1u(s)\phi(\mbox{\boldmath $\Phi$}^{-1}(s))\,ds,
\qquad u(0)=u(1)=0.$$
A general solution of this equation is
$$u(t)=c_0\int\limits_{\frac 12}^t\frac {\sin(\omega(\tau-t))\,d\tau}
{\omega\phi(\mbox{\boldmath $\Phi$}^{-1}(\tau))}+
c_1\cos(\omega t)+c_2\sin(\omega t),\qquad \omega=\mu^{1/2}.
$$
Substituting it into the boundary conditions and into (3.5) after some calculations
we derive the equation for eigenvalues\footnote {In \cite{KKW} this equation
is written in different but equivalent form.}:
$$\det\left[\begin{array}{ccl}
0 & \int\limits_0^1\mbox{\boldmath $\Phi$}^{-1}(\tau)\cos(\omega\tau)\,d\tau
& -\frac 1{\omega^2Q}+ \int\limits_0^1\int\limits_0^{\tau}
\mbox{\boldmath $\Phi$}^{-1}(\tau)\mbox{\boldmath $\Phi$}^{-1}(t)\frac
{\sin(\omega(t-\tau))}{\omega}\,dtd\tau\\
\sin(\omega) & 0 & \phantom {+\frac 1{Q^2}}-\frac 1{\omega}
\int\limits_0^1\mbox{\boldmath $\Phi$}^{-1}(\tau)\cos(\omega\tau)\,d\tau\\
\cos(\omega) & 1 & \phantom {+\frac 1{Q^2}-} \frac 1{\omega}
\int\limits_0^1\mbox{\boldmath $\Phi$}^{-1}(\tau)\sin(\omega\tau)\,d\tau
\end{array}\right]=0.$$
For $Q=0$ passing to the limit gives a natural result $\sin(\omega)=0$. We note
that in this example also a half of eigenfunctions for the original process (namely,
$\sin(2\pi kt)$, $k\in\mathbb N$) are orthogonal to $\psi$ in $L_2(0,1)$. Therefore,
they do not depend on $Q$.\medskip

Now we describe a case where the KL-expansion for the transformed process
can be constructed trivially. Let $\varphi=u_m$ be an eigenfunction of the covariance
$G_X$. Then $\psi=\lambda_m u_m$, where $\lambda_m$ is corresponding eigenvalue.
Hence all the eigenfunctions $u_k$, $k\ne m$, are orthogonal to $\psi$, and therefore,
$$\begin{array}{ll}
\widetilde\lambda_k=\lambda_k,\qquad&\widetilde u_k=u_k,\qquad k\ne m,\\
\widetilde\lambda_m=\lambda_m(1-q\alpha)^2,\qquad&\widetilde u_m=u_m.
\end{array}$$

Now we establish the relation which simplifies (2.8) for the processes under consideration.\medskip

{\bf Theorem 3}. Let the covariance $G_X$ be the Green function for an operator of the form
(3.1), and let $\widehat\alpha= \frac 1q$. If $\varphi\in L_2(0,1)$ then, as $\ep\to0$,
$${\mathbb P}\{\|{\cal X}_{{\bf \varphi}, \widehat\alpha}\|\leq\ep\} \sim
\frac {\sqrt {q}}{\|\varphi\|_2}\cdot
 \left(2\ell\sin(\textstyle\frac \pi{2\ell})\ep^2\right)^
{-\frac {\ell}{2\ell-1}}\cdot {\mathbb P}\{\|X\|\leq\ep\}.\eqno(3.12)
$$

{\bf Proof}. It is shown in \cite[Theorem 1.2]{Na1} that the process $X$ satisfies the
relation\footnote{The case of the operator $L_X$ with "separated"\ boundary conditions was
considered earlier in \cite[\S7]{NN}.}
$$F(r)={\mathbb P}\{\|X\|\leq\sqrt{r}\}\sim {\cal C}\cdot r^{\beta}
\exp\left(-{\mathfrak D} r^{-d}\right),\qquad r\to0,\eqno(3.13)$$
where $d=\frac 1{2\ell-1}$,
${\mathfrak D}=\frac 1{2d}\big(2\ell \sin(\frac {\pi}{2\ell})\big)^{-d-1}$
while the values of $\cal C$ and $\beta$ are now inessential for us.

The behavior of the distribution density $F'(r)$ for small $r$ was studied in
\cite[Theorem 3]{Lf2} in highly general situation. In our case, see also the proof
of Theorem 6.2 \cite{NN}) this result of \cite{Lf2} can be rewritten as follows:
$$F'(r)\sim {\cal C}{\mathfrak D}d\cdot r^{\beta-d-1}
\exp\left(-{\mathfrak D} r^{-d}\right),\qquad r\to0.\eqno(3.14)$$
This easily means that the asymptotics (3.13) is differentiable w.r.t. $r$.

Substituting (3.14) and (2.12) into (2.10) we obtain
$$\widetilde F(r) \sim \frac {{\cal C}{\mathfrak D}d}{\|\varphi\|_2}\cdot
\sqrt {\frac {2q}{\pi}}\cdot \int\limits_0^r \frac
{x^{\beta-d-1}}{\sqrt {r-x}}\,\exp\left(-{\mathfrak D}x^{-d}\right)\,dx. $$
Changing the variable $x=r(1-y)$ we get
$$\widetilde F(r) \sim \frac {{\cal C}{\mathfrak D}d}{\|\varphi\|_2}\cdot
\sqrt {\frac {2q}{\pi}}\cdot r^{\beta-d-\frac 12}
\exp\left(-{\mathfrak D}r^{-d}\right)
\int\limits_0^1 \frac
{(1-y)^{\beta-d-1}}{\sqrt {y}}
\exp\left(-\frac {\mathfrak D}{r^d}\bigl((1-y)^{-d}-1\bigr)\right)dy.$$
It is easily seen that for $y\ge r^{d/2}$ the integrand is exponentially
small. Therefore, one can integrate only over the interval $[0,r^{d/2}]$. In this interval we have
$(1-y)^{\beta-d-1}\sim1$ and $(1-y)^{-d}-1\sim yd$. Changing the variable $y=r^dz$ we arrive at
\begin{multline*}
\widetilde F(r) \sim \frac {{\cal C}{\mathfrak D}d}{\|\varphi\|_2}\cdot
\sqrt {\frac {2q}{\pi}}\cdot r^{\beta-\frac {d+1}2}
\exp\left(-{\mathfrak D}r^{-d}\right)\int\limits_0^{r^{-d/2}}
\frac 1{\sqrt {z}}\exp\left(-d{\mathfrak D}z\right)dz\sim \\
\sim\frac {{\cal C}}{\|\varphi\|_2}\cdot
\sqrt {2q{\mathfrak D}d}\cdot r^{\beta-\frac{d+1}2}
\exp\left(-{\mathfrak D}r^{-d}\right)\sim
\frac {\sqrt {2q{\mathfrak D}d}}{\|\varphi\|_2}
\cdot r^{-\frac{d+1}2}\cdot F(r).
\end{multline*}
This gives (3.12).\hfill$\square$\medskip

In the examples considered earlier the small ball behavior of the original
processes is well known (the processes $W$ and $B$ are classical ones while the process
$\overline W_1$ was studied in \cite{BNO}). Applying Theorems 1 and 3 we arrive at\medskip

{\bf Proposition 1}. We have, as $\ep\to0$,
$$\begin{array}{rclll}
{\mathbb P}\{\|{\cal W}_{{\bf 1}, \alpha}\|\leq\ep\} & \sim &
\left[
\begin{array}{l}
 \frac {4\ep}{\sqrt{\pi}|1-\frac {\alpha}3|} \cdot
\exp\left(-{\textstyle \frac 18}\ep^{-2}\right),\\
\\
 \frac {2\ep^{-1}}{\sqrt{3\pi}}\cdot
\exp\left(-{\textstyle \frac 18}\ep^{-2}\right),
\end{array}
\right.
& \begin{array}{l}\alpha\ne3,\\ \\  \alpha=3;
\end{array}
\\
\\
{\mathbb P}\{\|{\cal B}_{{\bf 1}, \alpha}\|\leq\ep\} & \sim &
\left[
\begin{array}{l}
 \frac {2\sqrt{2}}{\sqrt{\pi}|1-\frac {\alpha}{12}|} \cdot
\exp\left(-{\textstyle \frac 18}\ep^{-2}\right),\\
\\
\frac {\ep^{-2}}{\sqrt{6\pi}}\cdot
\exp\left(-{\textstyle \frac 18}\ep^{-2}\right),
\end{array}
\right.
&\begin{array}{l} \alpha\ne12,\\ \\ \alpha=12;
\end{array}
\\
\\
{\mathbb P}\{\|{\cal B}_{\varphi, \alpha}\|\leq\ep\} & \sim &
\left[
\begin{array}{l}
 \frac {2\sqrt{2}}{\sqrt{\pi}|1-\frac {17\alpha}{5040}|} \cdot
\exp\left(-{\textstyle \frac 18}\ep^{-2}\right),
\\
\\
 \frac {\sqrt{17}\ep^{-2}}{2\sqrt{21\pi}}\cdot
\exp\left(-{\textstyle \frac 18}\ep^{-2}\right),
\end{array}
\right.
&\begin{array}{l}\alpha\ne\frac {5040}{17},\\ \\ \alpha=\frac {5040}{17},
\end{array}& \varphi(t)=t(1-t);
\\
\\
{\mathbb P}\{\|{\cal X}_{{\bf 1}, \alpha}\|\leq\ep\} & \sim &
\left[
\begin{array}{l}
\frac {4\sqrt{2}\ep^{-\frac 13}}{\sqrt {3\pi}|1-\frac {\alpha}{120}|}
\cdot\exp\left(-{\textstyle \frac 38}\ep^{-\frac 23}\right),\\
 \\
\frac {\ep^{-\frac 53}}{3\sqrt {5\pi}}
\cdot\exp\left(-{\textstyle \frac 38}\ep^{-\frac 23}\right),
\end{array}
\right.&
\begin{array}{l}\alpha\ne120,\\ \\ \alpha=120,
\end{array} & X(t)={\overline W}_1(t);
\\
\\
{\mathbb P}\{\|{\cal B}_{\varphi, \alpha}\|\leq\ep\} & \sim &
\begin{array}{l}\frac {2\sqrt{2}}{\sqrt{\pi}|1-\alpha|} \cdot
\exp\left(-{\textstyle \frac 18}\ep^{-2}\right),
\end{array} &
\begin{array}{l}\alpha\ne1,\end{array}
& \varphi(t)=\frac 1{\phi(\mbox{\boldmath $\Phi$}^{-1}(t))}.
\end{array}
$$

{\bf Remark}. In the example {\bf5}, we cannot apply Theorem 2 (and consequently, also
Theorem 3) for $\alpha=1$.

\section{Some generalizations}

The construction (1.4) can be generalized to a special class of distributions $\varphi$.

Let $\varphi\in D'({\cal O})$ satisfy $q\equiv{\mathbb E}|\big<\varphi,X\big>|^2<\infty$.
We define a set of Gaussian functions by the formula similar to (1.4):
$${\cal X}_{\varphi,\alpha}(x)=X(x)-\alpha \psi(x)\big<\varphi,X\big>,\quad
x\in\overline{\cal O},\eqno(4.1)$$
where $\psi(x)={\mathbb E}X(x)\big<\varphi,X\big>$.\medskip

{\bf Remark}. In terms of the operator theory, the restriction on $\varphi$ means
$\varphi\in (Im({\mathfrak G}^{1/2}))'$. This implies $\psi\in Im({\mathfrak G}^{1/2})$, and
$q=\big<\varphi,\psi\big>$. In the random processes theory, $Im({\mathfrak G}^{1/2})$ is called
the kernel of distribution for the process $X$. $\varphi$ is called a linear measurable functional
of $X$, see \cite[\S9]{Lf1}.\medskip

For the processes (4.1), the Main Lemma and Corollary 1 hold true. Instead of Corollary 2,
we have the following analogue.\medskip

{\bf Corollary 2'}. For $\widehat\alpha=\frac 1q$ the process
${\cal X}_{\varphi,\widehat\alpha}(x)$ and r.v. $\big<\varphi,X\big>$ are independent. Moreover,
$\big<\varphi, {\cal X}_{\varphi,\widehat\alpha}\big>=0$ a.s.\medskip

Further, under assumption (2.1) Theorem 1 is valid. If, in addition, the covariance
$G_X(t,s)$ satisfies the assumptions of \S3, then the algorithm for the construction of
KL-expansion for the process ${\cal X}_{\varphi,\alpha}$ also runs.\medskip

Now we give some examples.\medskip

{\bf Example 6}. Let $X=W$, and let $\varphi(t)=\delta(t-1)$. Then $\psi(t)=G_W(t,1)=t$,
$q=G_W(1,1)=1$, and
$${\cal G}_{\varphi,\alpha}(t,s)=\min\{t,s\}+Qts\qquad (Q=\alpha^2-2\alpha).$$
Thus, for $\alpha\in]0,2[$ the process ${\cal W}_{\varphi,\alpha}$
has the same distribution as the Brownian bridge from zero to zero
with the length $-\frac 1Q$, see \cite[4.4.20]{BoS}.

Theorem 1 gives for $\alpha\ne1$
$${\mathbb P}\{\|{\cal W}_{\varphi, \alpha}\|\leq\ep\}\sim \frac {4\ep}
{\sqrt{\pi}|1-\alpha|} \cdot \exp\left(-{\textstyle \frac 18}\ep^{-2}\right),
\qquad \ep\to0.$$
This result coincides with a particular case of \cite[Proposition 1.9]{Na} and
\cite[Theorem 4.1]{NP}.\medskip

{\bf Remark}. For $\alpha=1$ the process ${\cal W}_{\varphi,\alpha}$ has the same
distribution as the standard Brownian bridge. It is easy to see that the infinite
product in (2.9) diverges, and Theorem 2 is not applicable. The same is true for
the forthcoming examples.\medskip

{\bf Example 7}. Let us consider the integrated Wiener process
$$X(t)=W^{[0]}_1(t)=\int\limits_0^tW(s)\,ds.$$

If $\varphi(t)=\delta'(t-1)$ then $\psi(t)=-(G_{W_1})_s(t,1)=-\frac {t^2}2$, \
$q=(G_{W_1})_{st}(1,1)=1$, and
$${\cal G}_{\varphi,\alpha}(t,s)=G_{W_1}(t,s)+Q\frac {t^2s^2}4\qquad
(Q=\alpha^2-2\alpha).$$
Thus, the process ${\cal X}_{\varphi,\alpha}$ has the same distribution as
the integrated process from the example {\bf6}. Our Theorem 1 corresponds to the case $m=1$ in
\cite[Proposition 1.9]{Na}. The small ball behavior in the case $\alpha=1$ (the integrated
Brownian bridge) was considered in \cite[Proposition 1.6]{Na}.\medskip

Now let $\varphi(t)=\delta(t-1)$. Then $\psi(t)=G_{W_1}(t,1)=\frac {t^2}2-\frac {t^3}6$, and
$q=G_{W_1}(1,1)=\frac 13$. Theorem 1 gives for $\alpha\ne3$
$${\mathbb P}\{\|{\cal X}_{\varphi, \alpha}\|\leq\ep\}\sim \frac
{8\sqrt{6}\,\ep^{1/3}}{\sqrt{\pi}|3-\alpha|} \cdot
\exp\left(-{\textstyle \frac 38}\ep^{-2/3}\right),
\qquad \ep\to0.$$
For $\alpha=3$ the direct calculation shows that
${\cal G}_{\varphi,\alpha}(t,s)$ is the Green function of the boundary value problem
$$u^{\scriptscriptstyle IV}=\mu u;\qquad u(0)=u'(0)=u(1)=u''(1)=0.$$
Applying \cite[Theorem 1.4]{Na} we obtain
$${\mathbb P}\{\|{\cal X}_{\varphi,\alpha}\|\leq\ep\}\sim \frac
{4\sqrt{2}\,\ep^{-2/3}}{3\sqrt{\pi}} \cdot
\exp\left(-{\textstyle \frac 38}\ep^{-2/3}\right),
\qquad \ep\to0.$$

{\bf Example 8}. Now we consider the Slepian process \cite{Sl}
that is a stationary zero mean-value Gaussian process with the
covariance $G_S(t,s)= 1- |t-s|$, $t,s\in [0,1]$. Let
$\varphi(t)=\delta(t)+\delta(t-1)$. Then $\psi\equiv1$, $q=2$, and
$${\cal G}_{\varphi,\alpha}(t,s)=1+Q-|t-s|\qquad (Q=2\alpha^2-2\alpha).$$
Thus, the process ${\cal S}_{\varphi,\alpha}$ has the same
distribution as a generalized Slepian process $S^{(c)}$, $c=1+Q$
(see \cite{GL} and \cite[\S2]{Na1}). Note that for $Q\ge0$ also
the equality
$$S^{(c)}(t)\stackrel {d}{=}W(t+c)-W(t),\qquad 0\le t\le1$$
holds true.

The statement of Theorem 1 corresponds \cite[Theorem 2.1, part
2]{Na1}. The small ball behavior in the case $\alpha=1/2$ (i.e.
$c=1/2$) is considered in \cite[Theorem 2.1, part 1]{Na1}.\medskip

{\bf Example 9}. In a similar way one can show that the statement
of \cite[Proposition 1.9]{Na} for any $m\in\mathbb N$,
\cite[Theorem 2.2, part 2]{Na1} and some theorems of \cite{NP} can
be considered as particular cases of our Theorem 1. Thus, this
theorem provides a unified approach to many formulas obtained
earlier.\medskip

One can also consider multiparametric analogues of the transform
(1.4). We restrict ourselves to a simplest case.\medskip

Let the functions $\varphi_1, \varphi_2\in L_{1,loc}({\cal O})$
satisfy the assumption (1.3) and the "orthogonality"\
assumption\footnote {For the distributions $\varphi_1,\varphi_2\in
(Im({\mathfrak G}^{1/2}))'$ this assumption has the form
$\big<\varphi_1,{\mathfrak G}\varphi_2\big>=0$.}
$$\int\limits_{\cal O}\int\limits_{\cal O}
G_X(u,v)\varphi_1(u)\varphi_2(v)\,dudv=0.$$
Let us consider a set
of Gaussian functions
$${\cal X}_{\mbox {\boldmath$\varphi,\alpha$}}(x)=X(x)-
\alpha_1\psi_1(x)\int\limits_{\cal O} X(u)\varphi_1(u)\,du-
\alpha_2\psi_2(x)\int\limits_{\cal O} X(u)\varphi_2(u)\,du,\quad
x\in\overline{\cal O},\eqno(4.2)$$
where $\psi_k={\mathfrak
G}\varphi_k$, $k=1,2$.

It is easy to see that the function (4.2) has the covariance
$${\cal G}_{\mbox {\boldmath$\varphi,\alpha$}}(x,y)=G_X(x,y)+
Q_1\psi_1(x)\psi_1(y)+Q_2\psi_2(x)\psi_2(y),$$
where
$Q_k=q_k\alpha_k^2-2\alpha_k$, $q_k=\big<\varphi_k,\psi_k\big>$.
Therefore, one can obtain the $L_2$-small ball asymptotics of the
function ${\cal X}_{\mbox {\boldmath$\varphi,\alpha$}}$ for
$\alpha_k\ne\frac 1{q_k}$, $k=1,2$, applying Theorem 1 twice. In
the same way we can obtain analogues of other statements.

\section{Appendix}

{\bf Lemma 5.1}. Let us consider two sequences: $\mu_k>0$ and
$b_k\ge0$, $k\in\mathbb N$. Let $\sum_k\mu_k^{-1}<\infty$ and
$\sum_k b_k<\infty$. Then for any $\gamma_1,\gamma_2\in\mathbb R$,
as $R\to\infty$,
$$\frac 1{2\pi}\int\limits_0^{2\pi}\ln\biggl|\gamma_1+\gamma_2
\sum_{k=1}^{\infty}\,\frac {b_k}{1-\frac {\mu_k}{R\exp(i\theta)}}\,\biggr|\,
d\theta\ \longrightarrow\ \ln
\biggl|\gamma_1+\gamma_2\sum_{k=1}^{\infty}b_k\biggr|. \eqno(5.1)$$

{\bf Proof}. Given $\theta\in]0,2\pi[$, the expression $1-\frac
{\mu_k}{R\exp(i\theta)}$ is bounded away from zero. By the
Lebesgue Dominated Convergence Theorem, we can pass to the limit
under the sum sign. Thus, the integrand in (5.1) converges to
$\ln|\gamma_1+\gamma_2\sum_kb_k|$ for $\theta\in]0,2\pi[$, and the
convergence is uniform over any segment.

Further, the expression under the absolute value sign has only
simple zeros and poles. Therefore, the integrand in (5.1) has only
logarithmic singularities, and one can easily construct a summable
majorant. Applying again the Lebesgue Theorem, we complete the
proof.\hfill$\square$\medskip

\bigskip

I am grateful to Prof. M.A.~Lifshits and to Prof. Ya.Yu.~Nikitin
for stimulating discussions and for some references. Also I am thankful to
Prof. M.S.~Birman who pointed me out the paper \cite{Ba}.

\end{document}